\UseRawInputEncoding 

\documentclass{article}

\usepackage{proof}
\usepackage{latexsym}
\usepackage{amsfonts}
\usepackage{amssymb}
\usepackage{cite}

\newcommand{\brem}{\begin{remark}}
\newcommand{\erem}{\end{remark}}
\newcommand{\blem}{\begin{lemma}}
\newcommand{\elem}{\end{lemma}}
\newcommand{\bth}{\begin{theorem}}
\newcommand{\ethm}{\end{theorem}}
\newcommand{\benu}{\begin{enumerate}}
\newcommand{\eenu}{\end{enumerate}}
\newcommand{\bdes}{\begin{description}}
\newcommand{\edes}{\end{description}}
\newcommand{\bdf}{\begin{definition}}
\newcommand{\edf}{\end{definition}}
\newcommand{\bcor}{\begin{cor}}
\newcommand{\ecor}{\end{cor}}
\newcommand{\bprp}{\begin{proposition}}
\newcommand{\eprp}{\end{proposition}}
\newcommand{\bmlem}{\begin{mlemma}}
\newcommand{\emlem}{\end{mlemma}}
\newcommand{\bclm}{\begin{claim}}
\newcommand{\eclm}{\end{claim}}
\newcommand{\bprf}{{\bf Proof}.\hspace{2mm}}

\newcommand{\eprf}{\hspace*{\fill} $\Box$}

\newcommand{\beqn}{\begin{equation}}
\newcommand{\eeqn}{\end{equation}}
\newcommand{\beqnarr}{\begin{eqnarray}}
\newcommand{\eeqnarr}{\end{eqnarray}}
\newcommand{\beqnarrs}{\begin{eqnarray*}}
\newcommand{\eeqnarrs}{\end{eqnarray*}}

\newcommand{\spand}{\,\&\,}

\newtheorem{theorem}{Theorem}[section]
\newtheorem{definition}[theorem]{Definition}
\newtheorem{proposition}[theorem]{Proposition}
\newtheorem{lemma}[theorem]{Lemma}
\newtheorem{cor}[theorem]{Corollary}
\newtheorem{remark}[theorem]{Remark}
\newtheorem{mlemma}[theorem]{Main Lemma}
\newtheorem{claim}[theorem]{Claim}

\newcommand{\alp}{\alpha}

\newcommand{\veps}{\varepsilon}

\newcommand{\ome}{\omega}

\newcommand{\bet}{\beta}
\newcommand{\gam}{\gamma}
\newcommand{\Gam}{\Gamma}

\newcommand{\sig}{\sigma}
\newcommand{\Sig}{\Sigma}

\newcommand{\fal}{\forall}
\newcommand{\exi}{\exists}

\newcommand{\Rarw }{\Rightarrow}

\newcommand{\Lrarw}{\Leftrightarrow}

\newcommand{\la}{\langle}
\newcommand{\ra}{\rangle}

\title{
Provably well-founded strict partial orders
}
\author{Toshiyasu Arai
\\
Graduate School of Mathematical Sciences
\\
University of Tokyo
\\
3-8-1 Komaba, Meguro-ku,
Tokyo 153-8914, JAPAN
\\
tosarai@ms.u-tokyo.ac.jp
}
\date{}
\begin{document}
\maketitle
\begin{abstract}
In this note we show through infinitary derivations that
each provably well-founded strict partial order in ${\rm ACA}_{0}$ admits an embedding to
an ordinal$<\veps_{0}$.
\end{abstract}

\section{Provably well-founded relations}

A \textit{strict partial order} $\prec$ is an irreflexive  $\fal n(n\not\prec n)$ and
transitive $\fal n,m,k(n\prec m\prec k\to n\prec k)$, relation on $\ome$.
Let `$\prec \mbox{ is a strict partial order}$' denotes the formula
$\fal n(n\not\prec n) \land \fal n,m,k(n\prec m\prec k\to n\prec k)$.
$<_{\veps_{0}}$ denotes a standard $\veps_{0}$-order, while
$<_{\ome}$ the usual order on $\ome$.

\bth\label{th:1}
Assume ${\rm ACA}_{0}\vdash {\rm TI}(\prec)$ for a primitive recursive relation $\prec$.
Then there exist an ordinal $\alp_{1}<\veps_{0}$ and a primitive recursive function $f$ such that
${\rm I}\Sig_{1}$ proves 
\[
\prec \mbox{ is a strict partial order } \to \fal n,m\left(n\prec m \to f(n)<_{\veps_{0}}f(m)<_{\veps_{0}}\alp_{1}\right)
.\]
\end{theorem}

Theorem \ref{th:1} is shown in \cite{Arai98a} by modifying Takeuti's proof in \cite{TRemark,PT2}
in terms of Gentzen's finitary proof\cite{GBew}.
In this note we show Theorem \ref{th:1} through infinitary derivations.

\bcor
Assume ${\rm ACA}_{0}\vdash {\rm TI}(\prec)$ for a primitive recursive relation $\prec$.
Then there exists an extension $\prec^{\prime}$ of $\prec$ such that
$\prec^{\prime}$ is primitive recursive, a well order, and
${\rm ACA}_{0}\vdash {\rm TI}(\prec^{\prime})$.
\ecor
\bprf
Let $n\prec^{\prime}m:\Lrarw f(n)<_{\veps_{0}}f(m) \lor\left(
f(n)=f(m) \land n<_{\ome}m\right)$.
\eprf

\section{Proof}
Assume for a primitive recursive relation $\prec$,
${\rm ACA}_{0}\vdash {\rm TI}(\prec)$.
In what follows argue in ${\rm I}\Sig_{1}$, and assume that $\prec$ is a strict partial order.
There exists an ordinal $\alp_{0}<\veps_{0}$ such that, cf.\cite{A2020}
\beqn\label{eq:E}
\fal n\left[ \vdash^{\alp_{0}}_{0}E(n)\right]
\eeqn
where $\vdash^{\alp}_{c}\Gam$ designates that `there exists a (primitive recursive) infinitary derivation of $\Gam$
with $\ome$-rule and the following inferences $(prg)$ and $(Rep)$
\[
\infer[(prg)]{\vdash^{\alp}_{c}\Gam}
{
\{
\vdash^{\bet}_{c}\Gam,E(m)
\}_{m\prec n}
}
\]
where $\bet<_{\veps_{0}}\alp$, $E$ is a fresh predicate symbol and $(E(n))\in\Gam$.
The subscript $0$ in $\vdash^{\alp_{0}}_{0}\Gam$ indicates that  a witnessed derivation is cut-free.
\[
\infer[(Rep)]{\vdash^{\alp}_{c}\Gam}
{\vdash^{\bet}_{c}\Gam}
\]
where $\bet<_{\veps_{0}}\alp$.

Formally we understand by (\ref{eq:E}) the following fact.
There exist a primitive recursive tree $T\subset{}^{<\ome}\ome$
and a primitive recursive function $H$ such that to each node $\sig\in T$, 
a five data $H(\sig)=(seq(\sig),ord(\sig), rul(\sig), crk(\sig), num(\sig))$ are
assigned by $H$.
Let $\Gam=seq(\sig)$, $\alp=ord(\sig)$, $c=crk(\sig)$ and $n=num(\sig)$.
Then $H(\sig)$ indicates that a sequent $\Gam$ is derived by a derivation in depth at most $\alp$
with cut rank $c$. $J=rul(\sig)$ is the last inference.
\[
\infer[(J)]{\sig \vdash^{\alp}_{c}\Gam}
{
\{
\sig_{i} \vdash^{\bet_{i}}_{c}\Gam_{i}
\}_{i\in I}
}
\]
has to be locally correct with respect to inferences $(\lor), (\land), (\exi), (\fal), (cut), (prg)$
and $(Rep)$,
and $\bet_{i}<_{\veps_{0}}\alp$ for each $i$.
Moreover when $J=rul(\sig)=(prg)$, $(E(n))\in\Gam=seq(\sig)$ with $n=num(\sig)$ is the main
formula of the $(prg)$. 
Then\footnote{$H(\la\, \ra)$ is arbitrary for the root $\la\,\ra$ of the tree.} $H(\la n\ra)=(\{E(n)\}, \alp_{0}, rul(\la n\ra),0)$ for each $n$.
Although $T$ is not assumed to be well-founded,
$rul(\sig)$ is either $(prg)$ or $(Rep)$ for each $\sig\in T$.
Therefore $seq(\sig)\subset\{E(n):n\in\ome\}$.
Let us assume that
\[
\infer[(prg)]{\sig \vdash^{\alp}_{c}\Gam}
{
\{
\sig*\la m\ra \vdash^{\bet}_{c}\Gam,E(m)
\}_{m\prec n}
}
\quad
\infer[(Rep)]{\sig \vdash^{\alp}_{c}\Gam}
{\sig*\la 0\ra \vdash^{\bet}_{c}\Gam}
\]

First we define nodes $\sig_{m}\in T$ by induction on $m$ as follows.
Let $\bet_{m}=ord(\sig_{m})$ and $\Gam_{m}=seq(\sig_{m})$ and $J_{m}=rul(\sig_{m})$.
Namely $\sig_{m} \vdash^{\bet_{m}}_{0}\Gam_{m}$.
It enjoys
\beqn\label{eq:min}
\fal n((E(n))\in\Gam_{m}\Rarw m\preceq n)
\eeqn
\textbf{Case 1}.
$\lnot\exi n<_{\ome}m(m\prec n)$:
Then let $\sig_{m}=\la m\ra$.
This means that $\bet_{m}=\alp_{0}$ and $\Gam_{m}=\{E(m)\}$.
\\
\textbf{Case 2}.
$\exi n<_{\ome}m(m\prec n)$:
Let $n_{0}<_{\ome}m$ be the $<_{\ome}$-least number such that $m\prec n_{0}$ and
$\bet_{n_{0}}=\min_{<_{\veps_{0}}}\{\bet_{n}: n<_{\ome}m, \, m\prec n\}$.
Consider the last inference 
$J_{n_{0}}=rul(\sig_{n_{0}})$ in the derivation of 
$\sig_{n_{0}} \vdash^{\bet_{n_{0}}}_{0}\Gam_{n_{0}}$.
\\
\textbf{Case 2.1}. The last inference $J_{n_{0}}$ is a $(prg)$:
\[
\infer[(prg)]{\sig_{n_{0}} \vdash^{\bet_{n_{0}}}_{0}\Gam_{n_{0}}}
{
\{
\sig_{n_{0}}*\la n\ra \vdash^{\bet}_{0}\Gam_{n_{0}}, E(n)
\}_{n\prec n_{1}}
}
\]
where $\bet<_{\veps_{0}}\bet_{n_{0}}$ and $(E(n_{1}))\in\Gam_{n_{0}}$ with
$n_{1}=num(\sig_{n_{0}})$.
We have $m\prec n_{0}\preceq n_{1}$ by (\ref{eq:min}).
Then let $\sig_{m}=\sig_{n_{0}}*\la m\ra$.
Let $\bet_{m}=\bet$ and $\Gam_{m}=\Gam_{n_{0}}\cup \{E(m)\}$.
If $(E(n))\in\Gam_{n_{0}}$, then $m\prec n_{0}\preceq n$ by (\ref{eq:min}).
Hence (\ref{eq:min}) is enjoyed for $\sig_{m}$ since $\prec$ is assumed to be transitive.
\\
\textbf{Case 2.2}. The last inference $J_{n_{0}}$ is a $(Red)$:
\[
\infer[(Rep)]{\sig_{n_{0}} \vdash^{\bet_{n_{0}}}_{0}\Gam_{n_{0}}
}
{
\sig_{n_{0}}*\la 0\ra \vdash^{\bet}_{0}\Gam_{n_{0}}
}
\]
where $\bet<\bet_{n_{0}}$.
Then let $\sig_{m}=\sig_{n_{0}}*\la 0\ra$.
This means $\bet_{m}=\bet$ and $\Gam_{m}=\Gam_{n_{0}}$.
Again (\ref{eq:min}) is enjoyed for $\sig_{m}$ by the transitivity of $\prec$.

\blem\label{lem:1}
$\fal m\fal n<_{\ome}m\left[
m\prec n \Rarw \bet_{m}<_{\veps_{0}}\bet_{n}
\right]
$.
\elem
\bprf
In \textbf{Case 2}, if $n<_{\ome}m$ and $m\prec n$,
then $\bet_{m}<_{\veps_{0}}\bet_{n_{0}}\leq_{\veps_{0}}\bet_{n}$.
\eprf
\\

Now let us define $\alp_{1}=\ome^{\alp_{0}}$ and $f$ as follows.

\[
f(n)=\max_{<_{\veps_{0}}}\{\ome^{\bet_{n_{0}}}\#\cdots\#\ome^{\bet_{n_{\ell-1}}}\#\ome^{\bet_{n_{\ell}}}:
 \fal i<\ell( n_{i}\prec n_{i+1} \spand n_{i}<_{\ome}n_{\ell}=n)\}
\]
where $\#$ denotes the natural sum.
Note that $n_{i}\neq n_{j}$ for $i<j\leq\ell$ 
since $\prec$ is assumed to be a strict partial order.
The following Lemma \ref{lem:2} shows Theorem \ref{th:1}.

\blem\label{lem:2}
$\fal n,m\left[n\prec m \Rarw f(n)<_{\veps_{0}}f(m)<\ome^{\alp_{0}+1}=\alp_{1}\right]$.
\elem
\bprf
Let 
$n_{0},\ldots,n_{\ell-1}<_{\ome}n_{\ell}=n\prec m$ be such that
$n_{0}\prec\cdots\prec n_{\ell-1}\prec n_{\ell}$ and
\[
f(n)=\ome^{\bet_{n_{0}}}\#\cdots\#\ome^{\bet_{n_{\ell-1}}}\#\ome^{\bet_{n_{\ell}}}
.\]
Then $n_{i}\prec m$ and $n_{i}\neq m$.
Let $A=\{i\leq \ell: m<_{\ome}n_{i}\}$ and $B=\{i\leq\ell: n_{i}<_{\ome}m\}$.
Then $A\cup B=\{0,\ldots,\ell\}$ and $A\cap B=\emptyset$.
By Lemma \ref{lem:1} we obtain
$\fal i\in A(\bet_{n_{i}}<_{\veps_{0}}\bet_{m})$, and hence
\beqn\label{eq:A}
\sum\{\ome^{\bet_{n_{i}}} : i\in A\}<_{\veps_{0}}\ome^{\bet_{m}}
\eeqn
where
$\sum\{\alp_{0},\ldots,\alp_{n}\}=\alp_{0}\#\cdots\#\alp_{n}$.
On the other side let
\[
\gam:= \max_{<_{\veps_{0}}}\{\ome^{\bet_{m_{0}}}\#\cdots\#\ome^{\bet_{m_{k-1}}}:
\fal i<k(m_{i}\prec m_{i+1}\spand m_{i}<_{\ome}m_{k}=m)\}
\]
and $B=\{n_{i_{0}}\prec\cdots\prec n_{i_{\ell-1}}\}$.
Then $n_{i_{0}}\prec\cdots\prec n_{i_{\ell-1}}\prec m$ and $n_{j}<_{\ome}m$
for each $n_{j}\in B$ since $\prec$ is assumed to be transitive.
Therefore
\beqn\label{eq:B}
\sum\{\ome^{\bet_{n_{i}}} : i\in B\} \leq_{\veps_{0}} \gam
\eeqn
By (\ref{eq:B}) and (\ref{eq:A}) we conclude 
\[
f(n)=\sum\{\ome^{\bet_{n_{i}}} : i\in B\} \# \sum\{\ome^{\bet_{n_{i}}} : i\in A\} 
<_{\veps_{0}}\gam\#\ome^{\bet_{m}}=f(m)
.\]
\eprf
\\

When $\prec$ is elementary recursive, then so is $f$.
For almost all theories $T$, Theorem \ref{th:1} holds if the ordinal $\veps_{0}$ is replaced by
the proof-theoretic ordinal of $T$
provided that a reasonable ordinal analysis of $T$ is given.

\end{document}